\documentclass{birkjour}
\newcommand{\f}{\gamma}
\renewcommand{\phi}{\varphi}

\def\beginpf{\medskip\noindent {\bf Proof.} \;}
\def\remark{\medskip\noindent {\bf Remark.} \;}

\begin{document}

%-------------------------------------------------------------------------
% editorial commands: to be inserted by the editorial office
%
%\firstpage{1} \volume{228} \Copyrightyear{2004} \DOI{003-0001}
%
%
%\seriesextra{Just an add-on}
%\seriesextraline{This is the Concrete Title of this Book\br H.E. R and S.T.C. W, Eds.}
%
% for journals:
%
%\firstpage{1}
%\issuenumber{1}
%\Volumeandyear{1 (2004)}
%\Copyrightyear{2004}
%\DOI{003-xxxx-y}
%\Signet
%\commby{inhouse}
%\submitted{March 14, 2003}
%\received{March 16, 2000}
%\revised{June 1, 2000}
%\accepted{July 22, 2000}
%
%
%
%---------------------------------------------------------------------------
%Insert here the title, affiliations and abstract:
%

\title[Averaged wave operators]
{Averaged wave operators \\ and complex-symmetric operators}

%----------Author 1
\author[R.~Bessonov]{Roman Bessonov}
\address{%
Chebyshev Laboratory \br
Department of Mathematics and Mechanics \br
St.Petersburg State University\br
14th Line V.O. 29b, St.Petersburg 199178 \br
Russia \br 
{\it and}}  \br
\address{St.Petersburg Department of the Steklov Mathematical Institute \br
Fontanka 27, St.Petersburg 191023 \br
Russia}
\email{bessonov@pdmi.ras.ru}

\thanks{This work is partially supported by RFBR grant 14-01-00748. The
first author is also supported by JSC ``Gazprom Neft'' and by the Chebyshev
Laboratory (Department of Mathematics and Mechanics, St. Petersburg State
University) under RF Government grant 11.G34.31.0026.}

%----------Author 2
\author[V.~Kapustin]{Vladimir Kapustin}
\address{St.Petersburg Department of the Steklov Mathematical Institute \br
Fontanka 27, St.Petersburg 191023 \br
Russia}
\email{kapustin@pdmi.ras.ru}
%----------classification, keywords, date


\date{January 1, 2015}
%----------additions
% \dedicatory{}
%%% ----------------------------------------------------------------------

\begin{abstract}
We study the behaviour of sequences $U_2^n X U_1^{-n}$, where $U_1, U_2$ are unitary operators, whose spectral measures are singular with respect to the Lebesgue measure, and the commutator $XU_1-U_2X$ is small in a sense. The conjecture about the weak averaged convergence of the difference $U_2^n X U_1^{-n}-U_2^{-n} X U_1^n$ to the zero operator is discussed and its connection with complex-symmetric operators is established in a general situation. For a model case where $U_1=U_2$ is the unitary operator of multiplication by $z$ on $L^2(\mu)$, sufficient conditions for the convergence as in the Conjecture are given in terms of kernels of integral operators.
\end{abstract}

%%% ----------------------------------------------------------------------
\maketitle
%%% ----------------------------------------------------------------------
%\tableofcontents

In scattering theory (in the discrete-time scheme) the behaviour of the sequence of operators $U_2^n X U_1^{-n}$ is studied, where $U_1, U_2$ are unitary operators on Hilbert spaces, and $X$ is a so-called identification operator. The limits of the sequence as $n\to\pm\infty$, whenever they exist, are called the wave operators of past and future. It is often assumed that the commutator 
\[
K=XU_1-U_2X
\]
is small in a certain sense; if $U_1$ and $U_2$ act on the same space and $X=I$, the operator $U_2$ can be regarded as a perturbation of $U_1$. In the typical scheme of scattering theory, only the absolutely continuous part of the spectral measure (with respect to the Lebesgue measure) of the unitary operator $U_1$ is considered. Usually, the singular subspace is killed by the orthogonal projection $P_{ac}$ onto the absolutely continuous subspace of $U_1$, and the sequence $U_2^n X U_1^{-n}P_{ac}$ is considered instead of $U_2^n X U_1^{-n}$. Pearson's version of the Kato--Rosenblum theorem \cite{DP78} says (for the continuous-time scheme) that if the commutator $K$ belongs to the trace class, then the strong wave operators of past and future exist. For the discrete-time scheme this means that the sequence $U_2^n X U_1^{-n} P_{ac}$ has strong limits as $n\to\pm\infty$. 

Without the restriction to the absolutely continuous subspace the direct analog of this result fails, as one can see on the following example for the one-dimensional space. If $U_1\ne U_2$ and $X=I$, then $U_2^n X U_1^{-n}=\omega^n I$, where $|\omega|=1$, $\omega\ne 1$. The numbers $\omega^n$ run over the unit circle and thus the sequence has no limit. However, we have the averaged convergence. Given a sequence $x_0, x_1, x_2, \dots$, its Ces\`aro means are defined by the formula
\[
\frac{1}{N}\sum_{n=0}^{N-1}x_n.
\]
The Ces\`aro means of the sequence $\omega^n$,
\[
\frac{1}{N}\sum_{n=0}^{N-1}\omega^n=\frac{1}{N}\frac{1-\omega^N}{1-\omega},
\]
tend to 0 as $N\to+\infty$. This naturally leads us to the question about the averaged convergence in the general case. For every pair $U_1, U_2$ of unitary operators, whose difference has rank one, the averaged limits at $\pm\infty$ of the sequence $U_2^n X U_1^{-n}$ always exist in the weak operator topology, cf.~Theorem~7.1 from \cite{VK10}. On the contrary, for every unitary operator $U_1$ with singular continuous spectral measure there exists a unitary operator $U_2$, which is a rank-two perturbation of $U_1$ and such that the weak averaged convergence of $U_2^n U_1^{-n}$ fails \cite{VK12}. It was shown in \cite{RB11} that, although the averaged limits at~$\pm\infty$ may not exist, the Ces\`aro means of the difference $U_2^n X U_1^{-n}-U_2^{-n} X U_1^n$ tend to zero in the weak operator topology if the commutator $XU_1-U_2X$ has rank two. This result confirms a partial case of the following conjecture, cf.~Conjecture~1.2 from \cite{VK10}.

\medskip
\noindent {\bf Conjecture.} \ {\it
Assume that $U_1, U_2$ are unitary operators, whose spectral measures are singular with respect to the Lebesgue measure, and the commutator $K=XU_1-U_2X$ has finite rank. Then the Ces\`aro means of the sequence 
\begin{equation}\label{eq1}
U_2^n X U_1^{-n}-U_2^{-n} X U_1^n
\end{equation}
tend to the zero operator in the weak operator topology as $n\to +\infty$.
}

\medskip

In particular, the Conjecture contains the assertion that if one of the weak averaged wave operators of past or future exists, then the other exists as well, and they both coincide with their half-sum. This motivates the study of the sequences 
\begin{equation}\label{eq2}
U_2^n X U_1^{-n}+U_2^{-n} X U_1^n.
\end{equation}
In this article we examine the weak averaged convergence of the sequences \eqref{eq1} and~\eqref{eq2}.

\medskip

The general case of a pair of unitary operators $U_1, U_2$ formally reduces to a special case where $U_1=U_2$ is the operator of multiplication by the independent variable $z$ on a space $L^2(\mu)$, where $\mu$ is a Borel measure on the unit circle which is singular with respect to the Lebesgue measure; moreover, one can assume that $\mu$ has no point masses. Indeed, the strong convergence on the absolutely continuous spectral subspace of $U_1$ follows from the classical scattering theory; it is easy to check that on eigenvectors strong limits of the Ces\`aro means exist. Since we study the weak convergence, we can assume that the spectral multiplicities of $U_1, U_2$ equal 1, and by the spectral theorem we obtain a pair of Borel measures in the spectral representations of $U_1, U_2$. If the spectral measures of $U_1, U_2$ are mutually singular, then the Ces\`aro means of the sequence $U_2^n X U_1^{-n}$ tend to the zero operator as $n\to\pm\infty$ in the weak operator topology \cite{VK12a}. Thus, the question reduces to the case where the measures are mutually absolutely continuous, and one can think that they coincide; we get the operator of multiplication by $z$ on $L^2(\mu)$ for some Borel measure $\mu$ on the unit circle. Since the absolutely continuous, purely point and singular continuous parts of $\mu$ can be considered separately, it remains to consider the case where $\mu$ is singular continuous. For more details see \cite{VK12a}.

Notice that the averaged convergence of \eqref{eq1} and \eqref{eq2} depends only on the commutator $K$. Namely, if $X, \tilde X$ are operators such that 
\[
K=XU_1-U_2X=\tilde XU_1-U_2\tilde X,
\]
then the convergence for these two situations holds or fails simultaneously. Indeed, then the difference $X-\tilde X$ satisfies $(X-\tilde X)U_1=U_2(X-\tilde X)$, and the sequences $U_2^n X U_1^{-n}$ and $U_2^n \tilde X U_1^{-n}$ differ by a constant sequence:
\[
U_2^n X U_1^{-n}- U_2^n \tilde X U_1^{-n}=U_2^n (X-\tilde X) U_1^{-n}= X-\tilde X.
\]
Hence we can look for conditions for the averaged convergence in terms of the commutator~$K$.  

Apparently, the assumption about the finite rank of $K$ in the Conjecture is not typical for possible applications, and the more natural smallness assumption is that it belong to the Schatten--von Neumann class~$\mathfrak S_1$ or~$\mathfrak S_2$. We formulate the Conjecture for finite-rank commutators, because the question is open for rank three and higher. Nevertheless, here we work in the situation where $K \in \mathfrak S_2$ is a Hilbert--Schmidt operator, because in this case $K$ can be rewritten as an integral operator. Namely, an operator $K$ on $L^2(\mu)$ belongs to the class $\mathfrak S_2$ if and only if it coincides with an integral operator, whose kernel $k$ belongs to $L^2(\mu\times\mu)$:
\[
(Kf)(z)=\int k(\xi, z) f(\xi)\,d\mu(\xi), \qquad \int |k(\xi, z)|^2\,d\mu(\xi)\,d\mu(z)<\infty.
\]
Thus, it is natural to study the averaged convergence in terms of the kernel~$k$ of the commutator $K \in \mathfrak S_2$ viewed as an integral operator. 

\medskip

Given a Borel measure $\mu$, the space $L^0(\mu)$ is the set of measurable functions defined $\mu$-almost everywhere; functions that coincide at $\mu$-almost all points are viewed as the same function. Every function from $L^0(\mu)$ can be represented as a ratio of two functions from $L^\infty(\mu)$ with the function in the denominator that may vanish only on a set of zero $\mu$-measure. 

A finite Borel measure $\mu$ on the unit circle is called a Rajchman measure if the sequence of its Fourier coefficients $\hat\mu(n)$ tends to zero as $n\to\infty$. All Borel measures that are absolutely continuous with respect to the Lebesgue measure are Rajchman measures; but there also exist singular Rajchman measures. More information about Rajchman measures can be found, e.g., in~\cite{Lyons}.

\medskip

\noindent {\bf Theorem.} \ {\it
Let $\mu$ be a Borel measure on the unit circle having no point masses, and denote by $U$ the operator of multiplication by $z$ on $L^2(\mu)$. Suppose that $X$ is a bounded operator on $L^2(\mu)$, for which the commutator $K=XU-UX$ belongs to the Hilbert--Schmidt class $\mathfrak S_2$ and coincides with an integral operator with kernel $k\in L^2(\mu\times\mu)$. 

\smallskip

\noindent $1)$ \; If for some function $\f\in L^0(\mu)$, $\f\ne 0$ \ $\mu$-almost everywhere, $k$ satisfies
\begin{equation}\label{f01}
\f(z)k(\xi, z)=-\f(\xi)k(z,\xi) 
\end{equation}
for $\mu\times\mu$-almost all pairs $\xi, z$, then the Ces\`aro means of $U^n X U^{-n}-U^{-n} X U^n$ tend to the zero operator in the weak operator topology.

\smallskip

\noindent $2)$ \; If for some function $\f\in L^0(\mu)$, $\f\ne 0$ \  $\mu$-almost everywhere, $k$ satisfies
\begin{equation}\label{f01p}
\f(z)k(\xi, z)=\f(\xi)k(z,\xi) 
\end{equation}
for $\mu\times\mu$-almost all pairs $\xi, z$, then the sequence of the Ces\`aro means of $U^n X U^{-n}+U^{-n} X U^n$ has a limit in the weak operator topology.

\smallskip

If $\mu$ is a Rajchman measure, then the analogs of statements $1)$ and $2)$ are fulfilled for the sequences \eqref{eq1} and \eqref{eq2}, respectively, without averaging.
}

\medskip

In the special case where $\f\equiv 1$ we obtain antisymmetric and symmetric kernels, respectively. 

The important question appears about a description of the linear span of all kernels $k$ that correspond to commutators of the form $K=XU-UX$ and satisfy \eqref{f01}. If it coincides with the set of all commutators from the Hilbert--Schmidt class $\mathfrak S_2$, then the Conjecture will be true even for the class $\mathfrak S_2$. One can also ask the same question about the subclasses that are the intersections of the class under consideration with $\mathfrak S_1$ or with the set of all finite-rank operators. 

The linear span of the commutators that satisfy \eqref{f01p} cannot coincide with the set of kernels of all commutators because there exist rank-two commutators, for which the weak averaged wave operators fail to exist. In this case one can ask if the linear span coincides with the set of all commutators, for which the corresponding wave operators exist.

The last assertion of the Theorem about Rajchman measures leads us to a stronger version of the above Conjecture. Namely, if the spectral measure of $U_1$ is absolutely continuous with respect to a certain Rajchman measure on the circle, then one can expect that the sequence~\eqref{eq1} weakly tends to zero, that is, that we can get the convergence even without averaging.

The Rajchman measures are characterized \cite{Lyons} as finite Borel measures on the circle that vanish on all Weyl sets (see the definition below).  Therefore, the condition that the spectral measure of $U_1$ is absolutely continuous with respect to a Rajchman measure is equivalent to the fact that it vanishes on all Weyl sets.

Weyl sets are defined as follows. A Borel measure $\nu$ is said to be the asymptotic distribution of a sequence $\omega_n$ of points of the unit circle if $\nu$ is the $\ast$-weak limit of the measures
\[
\frac{1}{N}\sum_{n=1}^N \delta_{\omega_n},
\]
where $\delta_\omega$ stands for the Dirac measure at $\omega$. A Borel set $e$ on the unit circle is called a Weyl set if there exists an increasing sequence $(n_k)$ of positive integers such that for every $\omega\in e$ the sequence $\omega^{n_k}$ has an asymptotic distribution, which does not coincide with the Lebesgue measure on the circle.

\medskip

Now let $\mu$ be singular. Every commutator $K=(\cdot, \bar u_1)v_1-(\cdot, \bar u_2)v_2$ of rank two on $L^2(\mu)$ with nonzero functions $u_1, v_1$, $u_2, v_2$ satisfies condition~\eqref{f01}. Indeed, $K$ is the integral operator with kernel 
\[
k(\xi, z)=u_1(\xi)v_1(z)-u_2(\xi)v_2(z).
\] 
If $K=XU-UX$ for a bounded operator $X$, then by Theorem 6.1 from \cite{VK10} 
\[
u_1 v_1-u_2 v_2=0
\]
$\mu$-almost everywhere. Define 
\[
\f=\frac{u_1}{v_2}=\frac{u_2}{v_1}.
\]
We obtain
\[
\begin{aligned}
\f(z)k(\xi, z)=&\frac{u_2(z)}{v_1(z)}\,u_1(\xi)v_1(z)-\frac{u_1(z)}{v_2(z)}\,u_2(\xi)v_2(z) \\
=&u_2(z)u_1(\xi)-u_1(z)u_2(\xi) \\
=&\frac{u_1(\xi)}{v_2(\xi)}\,u_2(z)v_2(\xi)-\frac{u_2(\xi)}{v_1(\xi)}\,u_1(z)v_1(\xi) \\
=&-\f(\xi)k(z, \xi).
\end{aligned}
\]

This proof cannot be directly generalized to the case of rank-three commutators and this shows that the general case of finite-rank commutators is much more complicated than the case of rank two. Moreover, for rank-three commutators property \eqref{f01} fails: for kernels $k$ satisfying $k(\xi, z)=k(z, \xi)$, this property can be fulfilled only in trivial cases. For instance, \eqref{f01} fails for the kernel $k(\xi, z)=1-{\rm Re}\,(\bar\xi z)$. In this case
\[
K=(\cdot, 1)1-\frac{1}{2}\Big( (\cdot, \bar z)\bar z+(\cdot, z)z\Big),
\]
and $K=XU-UX$ for $X=\frac{1}{2}\big( (\cdot, z)1-(\cdot, 1)\bar z\big)$. Then $X$ is the sum of two rank-one operators, and, correspondingly, $K$ is the sum of two rank-two operators; each of them satisfies \eqref{f01}, but with different functions $\f$.

\bigskip

The proofs of our results will be based on a lemma about complex-symmetric operators. For a definition of this class of operators, a Hilbert space involution $C$ should be fixed. The typical example is the complex conjugation on an $L^2$-space, $Cf=\bar f$. But we want to cover a more general situation, where $C$ is a mapping from a Hilbert space to another one. Namely, let $C: H_1\to H_2$ be a real-linear mapping such that
\[
(Ch, C\tilde h)=(\tilde h, h)
\]
for any $h, \tilde h\in H_1$. We have $\|Ch\|=\|h\|$ and $C(ih)=-iCh$, hence the image $\mathcal C=\{Ch: \; h\in H_1\}$ of $C$ is a closed linear subspace of $H_2$. 

An operator is said to be $C$-symmetric if $CT^\ast C=T$. Under the assumption that $C$ is surjective, the mapping $T\mapsto CT^\ast C$ is an involution on the set of all operators from $H_1$ to $H_2$.  In the partial case where $H_1=H_2$ and $C^2=I$, we obtain the known defintion; see, e.g., \cite{Putinar}. 

Also consider a pair of unitary operators $U_1: H_1\to H_1$, $U_2: H_2\to H_2$. We assume that
\[
CU_1=U_2^{-1}C,
\]
whence 
\begin{equation}\label{nn}
CU_1^n=U_2^{-n} C
\end{equation}
for any integer $n$. Therefore, $\mathcal C$ is a reducing subspace of the unitary operator $U_2$. In fact, the restriction of $U_2$ to $\mathcal C$ is unitarily equivalent to $U_1$. Indeed, if $e\in H_1$, by \eqref{nn} $C$ maps the reducing subspace of $U_1$ generated by $e$ onto the reducing subspace of $U_2$ generated by $Ce$, and the mapping 
\begin{equation}\label{mmw}
\sum c_n U_1^n e\mapsto \sum c_n U_2^n Ce
\end{equation}
is an isometry:
\[
\begin{aligned}
(U_2^m Ce,  U_2^n Ce)=(Ce, U_2^{n-m} Ce)=&(Ce, CU_1^{m-n} e) \\
=&(U_1^{m-n} e, e)=(U_1^m e, U_1^n e).
\end{aligned}
\]
Thus, if $H_1$ is a direct sum of singly-generated subspaces reducing $U_1$, then $\mathcal C$ is the direct sum of the images of these subspaces, which reduce $U_2$; the claim follows.

\medskip

It will be convenient to work with formal sums $\sum_{m=p}^q$, where $p>q$; they are naturally defined by 
\[
\sum_{m=p}^q=\sum_{m=p}^N-\sum_{m=q+1}^{N}=-\sum_{m=q+1}^{p-1},
\]
where $N$ is an arbitrary sufficiently large positive integer.
In particular, we have
\[
\sum_{m=p}^{p-1}=0.
\]

Although here we apply the Lemma below only to the space $L^2(\mu)$, we formulate it for a pair of spaces, because this situation is closer to possible applications, while the proof is no more difficult.

\pagebreak

\noindent {\bf Lemma.} \; {\it
Let $C$ be as above. For a linear operator $X: H_1\to H_2$ define $K=XU_1-U_2X$. Take $e\in H_1$ and set $\bar e=Ce\in H_2$. Fix integers $k, l$.

$1)$ \; If $CK^\ast C=-K$, then
\begin{equation}\label{fl1}
\Big( (U_2^n XU_1^{-n}-U_2^{-n}XU_1^n)\,U_1^k e, U_2^l \bar e\Big)
=\sum_{m=1}^{k+l} (U_2^nKU_1^{-n} \, U_1^{k-m} e, U_2^{l-m+1}\bar e).
\end{equation}
%was:
%\[
%\Big( (U_2^{-n}XU_1^n-U_2^n XU_1^{-n})\,U_1^k e, U_2^l \bar e\Big)
%=\sum_{m=1}^{k+l} (U_2^{-n}KU_1^n \, U_1^{k-m} e, U_2^{l-m+1}\bar e).
%\]

$2)$ \; If $CK^\ast C=K$, then
\begin{equation}\label{fl2}
\begin{aligned}
\Big( (U_2^nXU_1^{-n}&+U_2^{-n} XU_1^n)\,U_1^k e, U_2^l\bar e\Big)= \\
=\Big((X+&U_2^{k+l}XU_1^{-(k+l)})\,U_1^k e, U_2^l\bar e\Big)+\sum_{m=1}^{k+l}(U_2^n K U_1^{-n}\, U_1^{k-m} e, U_2^{l-m+1}\bar e). 
\end{aligned}
\end{equation}
}

\medskip

Thus, in both cases the convergence of the left-hand side as $n\to\infty$ and convergence of the Ces\`aro means reduce to those for finite sums of the form $(U_2^nKU_1^{-n}\, h_1, h_2)$. 

\medskip
\noindent {\bf Corollary.} \; {\it
In the conditions of the above Lemma, assume that the linear spans of the vectors $U_1^ne$ and $U_2^n\bar e$, where $n$ runs over the set of all integers, are dense in $H_1$ and $H_2$, respectively.

If $CK^\ast C=-K$ and the sequence $U_2^nKU_1^{-n}$ (the sequence of its Ces\`aro means) has a limit (tends to the zero operator) in the weak operator topology, then the same is true for the sequence $U_2^nXU_1^{-n}-U_2^{-n} XU_1^n$.

If $CK^\ast C=K$ and the sequence $U_2^nKU_1^{-n}$ (the sequence of its Ces\`aro means) has a limit in the weak operator topology, then the same is true for the sequence $U_2^nXU_1^{-n}+U_2^{-n} XU_1^n$.
}

\beginpf  Relations \eqref{fl1} and \eqref{fl2} imply the required weak convergence for pairs of vectors from dense sets, and the statements follow from the uniform boundedness of the sequences. \qed

\medskip

Below we shall see that if $CK^\ast C=K$, for the weak convergence of the Ces\`aro means of the sequence $U_2^nXU_1^{-n}+U_2^{-n} XU_1^n$ the assumption about the averaged convergence of $U_2^nKU_1^{-n}$ is superfluous.

\medskip

\noindent{\bf Proof of the Lemma.} \; 
For every pair of integers $p, q$ we have
\begin{equation}\label{f16}
\begin{aligned}
U_2^pXU_1^{-p}-U_2^q XU_1^{-q}
=&\sum_{m=p}^{q-1} (U_2^{m} XU_1^{-m}-U_2^{m+1} XU_1^{-m-1})  \\
=&\sum_{m=p}^{q-1} U_2^{m}KU_1^{-m-1}.
\end{aligned}
\end{equation}
Hence for a fixed pair of integers $k, l$ we obtain
\begin{equation}\label{f17}
\begin{aligned}
\Big( (U_2^pXU_1^{-p}-U_2^q XU_1^{-q})\,U_1^k e, U_2^l &\bar e\Big)
=\sum_{m=p}^{q-1} (U_2^{m}KU_1^{-m-1}\,U_1^k e, U_2^l \bar e) \\
=&\sum_{m=p}^{q-1} (KU_1^{k-m-1} e, U_2^{l-m} \bar e)=\sum_{m=p}^{q-1}\eta_m,
\end{aligned}
\end{equation}
where
\[
\eta_m=(KU_1^{k-m-1}e, U_2^{l-m}\bar e).
\]
Rewrite relation~\eqref{fl1} in terms of $\eta_m$. By \eqref{f17},
\[
\Big( (U_2^{-n}XU_1^n-U_2^n XU_1^{-n})\,U_1^k e, U_2^l \bar e\Big)=\sum_{m=-n}^{n-1}\eta_m,
\]
and for the right-hand side of \eqref{fl1} we obtain
\[
\begin{aligned}
\sum_{m=1}^{k+l} (U_2^nKU_1^{-n} \, U_1^{k-m} e, U_2^{l-m+1}\bar e)
=\sum_{m=1}^{k+l} (K& U_1^{k-m-n} e, U_2^{l-m-n+1}\bar e) \\
=&\sum_{m=1}^{k+l} \eta_{m+n-1}=\sum_{m=n}^{n+k+l-1} \eta_m.
\end{aligned}
\]
Thus, \eqref{fl1} can be rewritten as 
\[
-\sum_{m=-n}^{n-1}\eta_m=\sum_{m=n}^{n+k+l-1} \eta_m,
\]
or, equivalently,
\begin{equation}\label{f1111}
\sum_{m=-n}^{n+k+l-1}\eta_m=0.
\end{equation}
We have
\[
\begin{aligned}
(KU_1^p e, U_2^q\bar e)=&(U_1^p e,  K^\ast U_2^q C e)=(C K^\ast U_2^q C e, CU_1^p e) \\
=&(C K^\ast C U_1^{-q}e, U_2^{-p}C e)=(C K^\ast C \, U_1^{-q}e, U_2^{-p}\bar e).
\end{aligned}
\]
If $CK^\ast C=-K$, then
\begin{equation}\label{f332}
\begin{aligned}
\eta_m=(KU_1^{k-m-1}e, U_2^{l-m}&\bar e)=(C K^\ast C \, U_1^{m-l}e, U_2^{m-k+1}\bar e) \\
=&-(K \, U_1^{m-l}e, U_2^{m-k+1}\bar e)=-\eta_{k+l-1-m},
\end{aligned}
\end{equation}
that is, $\eta_p=-\eta_q$ whenever $p+q=k+l-1$. Then \eqref{f1111} easily follows.

To prove the second assertion, rewrite \eqref{fl2} in the form
\[
\begin{aligned}
((U_2^{-n}XU_1^n-X))\,U_1^k e, U_2^l\bar e)-((U_2^{k+l}XU_1^{-(k+l)}-&U_2^n XU_1^{-n})\,U_1^k e, U_2^l\bar e)= \\
=\sum_{m=1}^{k+l}&(K U_1^{k-m-n} e, U_2^{l-m+1-n}\bar e). 
\end{aligned}
\]
By \eqref{f17} this is equivalent to
\begin{equation}\label{f1112}
\sum_{m=-n}^{-1}\eta_m-\sum_{m=k+l}^{n-1}\eta_m=\sum_{m=1}^{k+l}\eta_{m+n-1}.
\end{equation}
If $CK^\ast C=K$, then, similarly to \eqref{f332}, we have
\begin{equation}\label{f1113}
\eta_m=\eta_{k+l-1-m},
\end{equation}
whence
\[
\sum_{m=-n }^{-1}\eta_m=\sum_{m=k+l}^{k+l-1+n}\eta_m
=\sum_{m=k+l}^{n-1}\eta_m+\sum_{m=n}^{k+l-1+n}\eta_m.
\]
The second sum in the right-hand side coincides with the sum in the right-hand side of \eqref{f1112}. Thus, \eqref{f1112} is proved and the proof of the Lemma is complete.
\qed

\medskip

\noindent{\bf Proof of the Theorem.} \; 
At first, assume that $\f\equiv 1$ and consider the mapping $C$ on $L^2(\mu)$, $Cf=\bar f$; define $e\equiv 1$. The assumption that $K$ is the integral operator with kernel $k$ that satisfies $k(\xi, z)=-k(z, \xi)$ implies
\begin{equation}\label{f29}
\begin{aligned}
(CK^\ast C\, f)(z)=\overline{\int\overline{k(z, \xi)}\cdot\overline{f(\xi)}\,d\mu(\xi)} 
=\int k(z, \xi)\, f(\xi)&\,d\mu(\xi) \\
=-\int k(\xi, z)\, f(\xi)\,d\mu(\xi)=&-(Kf)(z), 
\end{aligned}
\end{equation}
that is, the assumption $CK^\ast C=-K$ of the Lemma is fulfilled. Similarly, the relation $k(\xi, z)=k(z, \xi)$ yields $CK^\ast C=K$.

\smallskip

The strong convergence of the sequence $U^nKU^{-n}$ in the case where $K$ is a compact operator and the spectral measure of $U$ is absolutely continuous is a well known fact. For instance, the continuous-time analog of it was used in \cite{DP78} even without a proof. Here we need a kindred result in a slightly more general form.

The set of all compact operators coincides with the uniform closure of the linear span of all rank-one operators. Since $K$ is a compact operator, it suffices to prove the convergence as in the Corollary to the Lemma with $K$ replaced by an arbitrary rank-one operator $\tilde K=(\cdot, \bar u)v$. For $f\in L^2(\mu)$ we obtain
\[
(U^n\tilde KU^{-n}f)(z)=z^n v(z)\int \bar\xi^n f(\xi)u(\xi)\,d\mu(\xi)=z^n v(z)\cdot \widehat{fu\,\mu}(n),
\]
whence the norm of $U^n\tilde KU^{-n}f$ equals $\|v\|\cdot |\widehat{fu\,\mu}(n)|$. It is well known that if $\mu$ is a Rajchman measure and $w\in L^1(\mu)$, then the Fourier coefficients $\widehat{w\mu}(n)$ also tend to zero (this is obvious for trigonometric polynomials $w$, and they are dense in $L^1(\mu)$). Since $fu\in L^1(\mu)$, we obtain $\widehat{fu\,\mu}(n)\to 0$ for Rajchman measures $\mu$. For an arbitrary Borel measure $\mu$ having no point masses, the Ces\`aro means of the absolute values of $\widehat{fu\,\mu}(n)$ tend to zero by the Wiener theorem \cite{wiener}. Thus, we get the strong convergence of the sequence $U^n KU^{-n}$ for Rajchman measures, and the strong convergence of its Ces\`aro means for arbitrary continuous measures, which proves the Theorem in the case $\f\equiv 1$.

Now let $\f\not\equiv 1$. Take $\phi_1, \phi_2\in L^\infty(\mu)$ such that $\phi_1, \phi_2\ne 0$ \ $\mu$-almost everywhere. The proof of the Theorem is based upon two observations:

\smallskip

\noindent $a$) \ the averaged convergence as in the Theorem is equivalent to that for $\phi_2(U)X\phi_1(U)$ in place of $X$;

\smallskip

\noindent $b$) \ relations \eqref{f01} and \eqref{f01p} are equivalent to their modified versions, where $X$ is replaced by $\phi_2(U)X\phi_1(U)$, and $\f$ is replaced by $\frac{\phi_1}{\phi_2}\f$.

\smallskip

Assertion $a$) follows from the relation
\[
(U^n\,\phi_2(U)X\phi_1(U)\,U^{-n}h, \tilde h)=(U^nX U^{-n}\, \phi_1 h, \bar\phi_2 \tilde h).
\]
The direct implication is obvious, the converse follows from the convergence on dense sets, because the operators under consideration are uniformly bounded.

Since $K=XU-UX$, the commutator of $\phi_2(U)X\phi_1(U)$ with $U$ is the operator $\phi_2(U)K\phi_1(U)$. The function $k$ is the kernel of the integral operator $K$, hence $\phi_2(U)K\phi_1(U)$ is the integral operator with kernel $\phi_2(z)k(\xi, z)\phi_1(\xi)$. The fact that $k$ satisfies \eqref{f01} or \eqref{f01p} with parameter $\f\in L^0(\mu)$ is equivalent to
\[
\begin{aligned}
\frac{\phi_1(z)}{\phi_2(z)}\f(z)\cdot & \phi_2(z)k(\xi, z)\phi_1(\xi)
=\phi_1(z)\cdot \f(z)k(\xi, z)\cdot\phi_1(\xi) \\
=&\mp \, \phi_1(z)\cdot \f(\xi)k(z, \xi)\cdot\phi_1(\xi) 
=\mp\,\frac{\phi_1(\xi)}{\phi_2(\xi)}\f(\xi)\cdot\phi_2(\xi)k(z, \xi)\phi_1(z),
\end{aligned}
\]
and $b$) follows.

Write $\f=\frac{\phi_2}{\phi_1}$ with $\phi_1, \phi_2\in L^\infty(\mu)$. By $b$) the commutator of the operator $\phi_2(U)X\phi_1(U)$ and $U$ is the integral operator with kernel satisfying \eqref{f01} or \eqref{f01p} with $\f=1$. Now the Theorem follows from $a$). \qed

\medskip

The Theorem can be proved a little bit more directly in terms of the kernels of the integral operators. Then the calculations will remind the proof for the case of rank two from \cite{RB11}, while the main steps will repeat those in the proof above. 

\medskip
\noindent {\bf Proposition.} \; {\it
If $CX^\ast C=X$, then $CK^\ast C=-K$. If $U_1$ has multiplicity 1 and $C$ is surjective, then the converse is also true.
}

\medskip

\beginpf
If $CX^\ast C=X$, then
\[
\begin{aligned}
CK^\ast C=C(XU_1-&U_2X)^\ast C=CU_1^{-1}X^\ast C-CX^\ast U_2^{-1}C \\
=& U_2 C X^\ast C-CX^\ast C U_1=U_2X-XU_1=-K.
\end{aligned}
\]
To prove the converse, construct a bounded opeator $Y$ such that $CYC=Y^\ast$ and $YU_1-U_2Y=K$. Namely, the sequence 
\[
Y_n=X-\frac{1}{2}(U_2^nXU_1^{-n}+U_2^{-n}XU_1^n)
\]
is uniformly bounded. Hence from the sequence of its Ces\`aro means one can select a subsequence that converges in the weak operator topology. Denote its limit by $Y$. We have
\[
\begin{aligned}
\left\|\left(\frac{1}{N}\sum_{n=0}^{N-1}U_2^nXU_1^{-n}\right)U_1-U_2\left(\frac{1}{N}\sum_{n=0}^{N-1}U_2^nXU_1^{-n}\right)\right\|&= \\
=\frac{1}{N}\,\|XU_1-U_2^NXU_1^{-N+1}&\| \leq \frac{2\|X\|}{N}\to 0,
\end{aligned}
\]
whence $YU_1-U_2Y=XU_1-U_2X=K$. Next, by formula~(\ref{f16}) we obtain
\[
\begin{aligned}
C(X-U_2^nX&U_1^{-n})^\ast C=C\left(\sum_{m=0}^{n-1}U_2^m K U_1^{-m-1}\right)^\ast C \\
=&\sum_{m=0}^{n-1} CU_1^{m+1}K^\ast U_2^{-m}C=\sum_{m=0}^{n-1} U_2^{-m-1}CK^\ast CU_1^m \\
=&-\sum_{m=0}^{n-1} U_2^{-m-1}KU_1^m=X-U_2^{-n}X U_1^n,
\end{aligned}
\]
and, similarly,
\[
C(X-U_2^{-n}XU_1^n)^\ast C=X-U_2^nXU_1^{-n}.
\]
Thus,
\[
\begin{aligned}
2\,CY_n^\ast C=C(X-&U_2^nXU_1^{-n})^\ast C+C(X-U_2^{-n}XU_1^n)^\ast C \\
=&(X-U_2^{-n}XU_1^n)+(X-U_2^nXU_1^{-n})=2\,Y_n
\end{aligned}
\]
and $CY^\ast C=Y$.

For $Z=X-Y$ we have $ZU_1=U_2 Z$. Suppose that the multiplicity of $U_1$ equals 1. Take a vector $e$ such that the vectors $U_1^n e$ generate the entire space. By the spectral theorem there exists a Borel measure $\mu$ such that $U_1$ is unitarily equivalent to the operator $U$ of multiplication by $z$ on $L^2(\mu)$, and $e\in H_1$ corresponds to $1\in L^2(\mu)$. The mapping \eqref{mmw} realizes the unitary equivalence of $U_1$ and $U_2$; thus one can think that $U_1=U_2=U$ and $Cf=\bar f$. To prove the Proposition, we must establish the implication
\begin{equation}\label{impl}
ZU=UZ \quad \implies \quad CZ^\ast C=Z;
\end{equation}
then 
\[
CX^\ast C=C(Y^\ast+Z^\ast)C=Y+Z=X.
\]
Indeed, since the multiplicity of $U$ is equal to 1, $ZU=UZ$ if and only if $Z=\phi(U)$ for some $\phi\in L^\infty(\mu)$, and then for $f\in L^2(\mu)$
\[
CZ^\ast Cf=\overline{\bar\phi\bar f}=\phi f=Zf,
\]
as required.  \qed
 
\medskip

Without the assumption that the multiplicity of $U$ equals 1, the implication \eqref{impl} fails and the reverse implication in the Proposition is not true.

\medskip

Quite similarly, the property $CX^\ast C=-X$ yields $CK^\ast C=K$. For the reverse implication we can say only that the operator $Y$, which is the limit of a subsequence of the Ces\`aro means of the operators $Y_n$ defined above, satisfies $CY^\ast C=-Y$ and $YU_1-U_2Y=K$. However, in the case of multiplicity 1 the operator having these properties is unique. Indeed, for operators $Z$ on $L^2(\mu)$ such that $ZU=UZ$, by \eqref{impl} the property $CZ^\ast C=-Z$ yields $Z=0$.

This argument shows that if $CK^\ast C=K$, then the Ces\`aro means of the operators $Y_n$ weakly tend to $Y$. Indeed, otherwise we can choose two subsequences weakly converging to limit operators, whose difference $Z\ne 0$ satisfies $ZU=UZ$ and $CZ^\ast C=-Z$. Since this is impossible, we obtain an improved version of one of the assertions from the Corollary to the Lemma: 

\smallskip

\noindent {\it if the closed linear span of the vectors $U_1^ne$ coincides with $H_1$, the closed linear span of the vectors $U_2^n Ce$ coincides with $H_2$, and $CK^\ast C=K$, then the Ces\`aro means of the sequence $U_2^nXU_1^{-n}+U_2^{-n} XU_1^n$ have a limit in the weak opeator topology}.

\medskip

\remark  Take a unimodular function $\gamma\in L^2(\mu)$. The mapping $C_\gamma$ on $L^2(\mu)$, $C_\gamma f=\bar\gamma\bar f$, also satisfies the properties of the involution $C$ listed above. For the involution $Cf=\bar f$, the condition $CK^\ast C=-K$ can be rewritten as $C_\gamma (KM_\gamma)^\ast C_\gamma=-KM_\gamma$, where $M_\gamma$ is the operator of multiplication by $\gamma$. Thus, the condition for $K$ and $C$ is equivalent to that for $KM_\gamma$ and $C_\gamma$. In our scheme the class of commutators, for which the averaged limit of the corresponding sequence $U^nXU^{-n}-U^{-n}XU^n$ exists, is closed with respect to bordering by operators of multiplication. Hence the classes of operators that can be obtained from the sets of the commutators $K$ satisfying $CK^\ast C=-K$ or $C_\gamma K^\ast C_\gamma=-K$ coincide, and it suffices to consider only the involution $f\mapsto\bar f$.


\begin{thebibliography}{8}
\bibitem{RB11} R.V.~Bessonov, \textit{The past and future wave operators on the singular spectrum.} Zap. Nauchn. Semin. POMI \textbf{389} (2011), 5--20 (Russian), English transl.: J. Math. Sci. \textbf{182} (2012), no.5, 587--594.
\bibitem{Putinar} S.R.~Garcia, M.~Putinar, \textit{Complex symmetric operators and applications}, Trans. Amer. Math. Soc. \textbf{358} (2006), no.3, 1285--1315.
\bibitem{VK10} V.V.~Kapustin, \textit{On wave operators on the singular spectrum.} Zap. Nauchn. Semin. POMI \textbf{376} (2010), 48--63 (Russian), English transl.: J. Math. Sci. \textbf{172}~(2011), no.2, 207--214.
\bibitem{VK12} V.V.~Kapustin, \textit{Averaged wave operators on the singular spectrum.} Funkts. Analiz i Prilozhen. \textbf{46} (2012) no.2, 24--36 (Russian), English transl.: Funct. Anal. and Its Appl. \textbf{46} (2012), no.2, 100--109.
\bibitem{VK12a} V.V.~Kapustin, \textit{Cauchy-type integrals and singular measures}, Algebra i Analiz \textbf{24}~(2012), no.5, 72--93 (Russian), English transl.: St.Petersburg Math. J. \textbf{24}~(2013), no.5, 743--757.
\bibitem{Lyons} R.~Lyons, \textit{Fourier-Stieltjes coefficients and asymptotic distribution modulo 1}, Ann. Math.  \textbf{122} (1985), no.1, 155--170.
%\textit {Seventy years of Rajchman measures}, Kahane Special Volume, J. Fourier Analysis and Applic., 363--377, 1995.
\bibitem{DP78} D.B.~Pearson, \textit{A generalization of the Birman trace theorem.} J. Funct. Anal. \ \ \textbf{28} (1978), no.2, 182-186.
\bibitem{wiener} N.~Wiener, \textit{The quadratic variation of a function and its Fourier coefficients.} \ J. of Math. and Phys. \textbf{3} (1924), 72--94.
\end{thebibliography}
\end{document}